\documentclass[10pt]{amsart}
\input{epsf}
\usepackage{amssymb,latexsym}
\topskip
\newcommand{\resp}[1]{{\rm [}respectively; #1{\rm ]}}

\numberwithin{equation}{section}
\newtheorem{theorem}{Theorem}[section]
\newtheorem{proposition}[theorem]{Proposition}
\newtheorem{corollary}[theorem]{Corollary}

\newtheorem{remark}[theorem]{Remark}
\newtheorem{lemma}[theorem]{Lemma}

\newtheorem{example}[theorem]{Example}

\begin{document}

\pagenumbering{arabic}
\pagestyle{headings}
\def\sof{\hfill\rule{2mm}{2mm}}
\def\ls{\leq}
\def\gs{\geq}
\def\SS{\mathcal S}
\def\qq{{\bold q}}
\def\txx{\left(\frac1{2\sqrt{x}}\right)}
\def\aa{\overline{1}}
\def\ab{\overline{2}}
\def\ac{\overline{3}}
\def\ad{\overline{4}}
\def\ae{\overline{5}}
\def\ak{\overline{k}}
\def\an{\overline{n}}
\def\EE{\frak E}
\def\OO{\frak O}
\def\MM{\frak M}
\def\SS{\frak S}
\def\sign{\mbox{{\rm sign}}}
\def\F{\mathcal F}
\def\vr{\varnothing}
\def\mn{\mbox{-}}
\def\sx{\left(\frac{1}{2x}\right)}
\def\ttx{\left(\frac{1}{2\sqrt{x}}\right)}

\title{Restricted even permutations and Chebyshev polynomials}
\maketitle
\begin{center}Toufik Mansour \footnote{Research financed by EC's
IHRP Programme, within the Research Training Network "Algebraic
Combinatorics in Europe", grant HPRN-CT-2001-00272}
\end{center}
\begin{center}
Department of Mathematics, Chalmers University of Technology,
S-412~96 G\"oteborg, Sweden

{\tt toufik@math.chalmers.se}
\end{center}
\section*{Abstract}
We study generating functions for the number of even (odd)
permutations on $n$ letters avoiding $132$ and an arbitrary
permutation $\tau$ on $k$ letters, or containing $\tau$ exactly
once. In several interesting cases the generating function depends
only on $k$ and is expressed via Chebyshev polynomials of the
second kind.

\noindent{2000 Mathematics Subject Classification}: Primary 05A05,
05A15; Secondary 30B70, 42C05
\section{Introduction}
The aim of this paper is to give analogies of enumerative results
on certain classes of permutations characterized by
pattern-avoidance in the symmetric group $\SS_n$. In the set of
even (odd) permutations we identify classes of restricted even
(odd) permutations with enumerative properties analogous to
results on permutations. More precisely, we study generating
functions for the number of even (odd) permutations avoiding $132$
and avoiding (or containing exactly once) an arbitrary permutation
$\tau\in\SS_k$. Moreover we consider statistics of the increasing
pattern. In the remainder of this section we present a brief
account of earlier works which motivated our investigation, we
give the basic definitions used throughout the paper, and the
organization of this paper.

Let $[p]=\{1,\dots,p\}$ denote a totally ordered alphabet on $p$
letters, and let $\alpha=(\alpha_1,\dots,\alpha_m)\in [p_1]^m$,
$\beta=(\beta_1,\dots,\beta_m)\in [p_2]^m$. We say that $\alpha$
is {\it order-isomorphic\/} to $\beta$ if for all $1\leq i<j\leq
m$ one has $\alpha_i<\alpha_j$ if and only if $\beta_i<\beta_j$.
For two permutations $\pi\in \SS_n$ and $\tau\in \SS_k$, an {\it
occurrence\/} of $\tau$ in $\pi$ is a subsequence $1\leq
i_1<i_2<\dots<i_k\leq n$ such that $(\pi_{i_1},\dots,\pi_{i_k})$
is order-isomorphic to $\tau$; in such a context $\tau$ is usually
called the {\it pattern\/}. We say that $\pi$ {\it avoids\/}
$\tau$, or is $\tau$-{\it avoiding\/}, if there is no occurrence
of $\tau$ in $\pi$, and we say that $\pi$ {\em contains} $\tau$
exactly $r$ times if there is $r$ different occurrences of $\tau$
in $\pi$. For example, the permutation $598376412\in\SS_9$ avoids
$123$ and contains $1432$ exactly twice.

While the case of permutations avoiding a single pattern has
attracted much attention, the case of multiple pattern avoidance
remains less investigated. In particular, it is natural, as the
next step, to consider permutations avoiding pairs of patterns
$\tau_1$, $\tau_2$. This problem was solved completely for
$\tau_1,\tau_2\in\SS_3$ (see \cite{SS}), and for $\tau_1\in\SS_3$
and $\tau_2\in\SS_4$ (see \cite{W}). Several recent papers
\cite{CW,MV1,Kr,MV2,MV3,MV4,BCS} deal with the case
$\tau_1\in\SS_3$, $\tau_2\in\SS_k$ for various pairs
$\tau_1,\tau_2$. Another natural question is to study permutations
avoiding $\tau_1$ and containing $\tau_2$ exactly $t$ times. Such
a problem for certain $\tau_1,\tau_2\in\SS_3$ and $t=1$ was
investigated in \cite{R}, and for certain $\tau_1\in\SS_3$,
$\tau_2\in\SS_k$ in \cite{RWZ,MV1,Kr,MV2,MV3,MV4}. For example,
several authors \cite{RWZ,MV1,Kr,BCS} have shown that generating
functions for the number $132$-avoiding permutations in $\SS_n$
with respect to number of occurrences of the pattern $12\ldots k$
can be expressed as either continued fractions or Chebyshev
polynomials of the second kind.

{\em Chebyshev polynomials of the second kind\/} (in what follows
just Chebyshev polynomials) are defined by
$U_r(\cos\theta)=\frac{\sin(r+1)\theta}{\sin\theta}$ for $r\geq0$.
The Chebyshev polynomials satisfy the following recurrence
$U_r(t)=2tU_{r-1}(t)-U_{r-2}(t)$ for $r\geq2$ together with
$U_0(t)=1$ and $U_1(t)=2t$. Evidently, $U_r(x)$ is a polynomial of
degree $r$ in $x$ with integer coefficients. Chebyshev polynomials
were invented for the needs of approximation theory, but are also
widely used in various other branches of mathematics, including
algebra, combinatorics, and number theory (see \cite{Ri}).
Apparently, for the first time the relation between restricted
permutations and Chebyshev polynomials was discovered  by Chow and
West in \cite{CW}, and later by Mansour and Vainshtein
\cite{MV1,MV2,MV3,MV4} and Krattenthaler \cite{Kr}. These results
are related to a rational function
\begin{equation}
R_k(x)=\frac{U_{k-1}\ttx}{\sqrt{x}U_k\ttx} \label{drk}
\end{equation}
for all $k\geq 1$. For example, $R_1(x)=1$,
$R_2(x)=\frac{1}{1-x}$, and $R_3(x)=\frac{1-x}{1-2x}$. It is easy
to see that for any $k$, $R_k(x)$ is rational in $x$ and satisfies
the following equation (see \cite{MV1,MV3,MV4}) for $k\geq1$,
\begin{equation}
R_k(x)=\frac{1}{1-xR_{k-1}(x)}. \label{rk}
\end{equation}


Let $\pi\in\SS_n$. The number of {\em inversions} of $\pi$ is
given by $|\{(i,j): \pi_i>\pi_j,\ 1\leq i<j\leq n\}|$. The {\em
sign} of $\pi$, $\sign(\pi)$, is given by the number of inversions
of $\pi$ modulo $2$ (equals $1$ if the number inversions of $\pi$
is given by even number, otherwise equals $-1$). We say $\pi$ is
an {\em even permutation} \resp{{\em odd permutation}} if
$\sign(\pi)=1$ \resp{$\sign(\pi)=-1$}. We say $\pi$ is an
involution if $\pi=\pi^{-1}$. We denote the set of all even
\resp{odd} permutations in $\SS_n$ by $\EE_n$ \resp{$\OO_n$}. For
example, if $\pi=4132\in\SS_4$ then the number of inversions of
$\pi$ equals $4$, so $\sign(\pi)=1$ and $\pi\in\EE_4$.

The paper of Simion and Schmidt \cite{SS} generalized and
considered for many directions. Here we give two examples. The
first one is paper of Chow and West~\cite{CW}, which had dealt
with three cases of avoiding $132$ or $123$ and avoiding
$\tau\in\SS_k$. In particular, they found the generating function
for the number of permutations in $\SS_n(132,12\ldots k)$ which
given by $R_k(x)$. The second one is papers of Guibert and Mansour
\cite{GM1,GM2}, which had dealt with the cases of avoiding $132$
(or containing exactly once) and avoiding $\tau\in\SS_k$ (or
containing exactly once). More precisely, the paper \cite{GM1}
dealt with the case of the generating function for number of
involutions in $\SS_n$ avoiding $132$ (or containing exactly once)
and avoiding $\tau\in\SS_k$ (or containing exactly once). The
paper \cite{GM2} dealt with the case of the generating function
for number of even (odd) involutions in $\SS_n$ avoiding $132$ (or
containing exactly once) and avoiding $\tau\in\SS_k$ (or
containing exactly once).

\begin{theorem}{\rm(see \cite{CW,GM1,GM2})} For all $k\geq0$,

{\rm(i)} The generating function for the number of permutations in
$\SS_n$ avoiding both $132$ and $12\ldots k$ is given by
            $$R_k(x).$$

{\rm(ii)} The generating function for the number of involutions in
$\SS_n$ avoiding both $132$ and $12\ldots k$ is given by
    $$I_k(x)=\frac{1}{xU_k\sx}\sum\limits_{j=0}^{k-1}U_j\sx.$$

{\rm(iii)} The generating function for the number of even
involutions in $\SS_n$ avoiding both $132$ and $12\ldots k$ is
given by
$$\sum\limits_{j=0}^{k-1}\left(x^j\biggl(1+\frac{x^2}{2}(R_{k-1-j}(x^2)+R_{k-1-j}(-x^2))I_{k-j}(x)\biggr)\prod\limits_{i=k-j}^k
R_i(-x^2)\right).$$
\end{theorem}
The above theorem invites the following question: Find explicitly
the generating function for the number even (odd) permutations
avoiding both $132$ and $12\ldots k$ in terms of Chebyshev
polynomials? In this paper we give a complete answer for this
question (see Subsection~\ref{ans1}).

%
As a consequence of \cite{MV2,MV4}, we present a general approach
to the study of even (odd) permutations avoiding $132$ and
avoiding an arbitrary pattern $\tau$ of length $k$, or containing
$\tau$ exactly once. We derive all the previously known results
for this kind of problems, as well as many new results.

The paper is organized as follows. The case of even (odd)
permutations avoiding both $132$ and $\tau$ is treated in
Section~\ref{sec2}. We derive a simple recursion for the
corresponding generating functions for general $\tau$. This
recursion can be solved explicitly for several interesting cases,
including $12\dots k$, $(d+1)(d+2)\ldots k12\ldots d$, and {\it
odd-wedge\/} patterns defined below. In particularly, we prove the
generating function for the number of even (odd) permutations
avoiding both $132$ and $\tau\in S_k(132)$ is a rational function
for every nonempty pattern $\tau$. Observe that if $\tau$ itself
contains $132$, then any $132$-avoiding permutation avoids $\tau$
as well, so in what follows we always assume that $\tau\in
S_k(132)$. The case of permutations avoiding $132$ and containing
$\tau$ exactly once is treated in Section~\ref{sec3}. Here again
we start from a general recursion, and then solve it for several
particular cases. Finally, in Section~\ref{sec4} we describe
several directions to extend and to generalize the results of the
pervious sections.

Most of the explicit solutions obtained in the next sections
involve Chebyshev polynomials of the second kind.
\section{Avoiding an arbitrary pattern}\label{sec2}
Consider an arbitrary pattern $\tau=(\tau_1,\dots,\tau_k)\in
\SS_k(132)$. Recall that $\tau_i$ is said to be a {\it
right-to-left maximum\/} if $\tau_i>\tau_j$ for any $j>i$. Let
$m_0=k,m_1,\dots,m_r$ be the right-to-left maxima of $\tau$
written from left to right. Then $\tau$ can be represented as
$$\tau=(\tau^0,m_0,\tau^1,m_1,\dots,\tau^r,m_r),$$
where each of $\tau^i$ may be possibly empty, and all the entries
of $\tau^i$ are greater than $m_{i+1}$ and all the entries of
$\tau^{i+1}$. This representation is called the {\it canonical
decomposition\/} of $\tau$. Given the canonical decomposition, we
define the $i$th {\it prefix\/} of $\tau$ by
$\pi^i=(\tau^0,m_0,\dots,\tau^i,m_i)$ for $1\leq i\leq r$ and
$\pi^0=\tau^0$, $\pi^{-1}=\emptyset$. Besides, the $i$th {\it
suffix\/} of $\tau$ is defined by
$\sigma^i=(\tau^i,m_i,\dots,\tau^r,m_r)$ for $0\leq i\leq r$ and
$\sigma^{r+1}=\emptyset$. Strictly speaking, prefixes and suffices
themselves are not patterns, since they are not permutations
(except for $\pi^r=\sigma^0=\tau$). However, any prefix or suffix
is order-isomorphic to a unique permutation, and in what follows
we do not distinguish between a prefix (or suffix) and the
corresponding permutation.

The set of all $T$-avoiding even \resp{odd} permutations in
$\SS_n$ we denote by $\EE_n(T)$ \resp{$\OO_n(T)$}. Let
$e_{\tau}(n)$ \resp{$o_{\tau}(n)$} be the cardinality of the set
$\EE_n(132,\tau)$ \resp{$\OO_n(132,\tau)$}. The corresponding
generating function let us denote by $\EE_{\tau}(x)$
\resp{$\OO_{\tau}(x)$}, that is,
$$\EE_\tau(x)=\sum_{n\geq 0} e_{\tau}(n)x^n\quad\left[\mbox{respectively; }\OO_\tau(x)=\sum_{n\geq 0} o_{\tau}(n)x^n\right].$$
The generating function for the number of permutations avoiding
both $132$ and $\tau$ we denote by $\F_\tau(x)$. Clearly,
\begin{equation}
\F_\tau(x)=\EE_\tau(x)+\OO_\tau(x).\label{eq3}
\end{equation}
The following proposition which is the base of all the results in
this section.

\begin{proposition}\label{gen1}
Let $k\geq 1$ and $n\geq 1$, then
$$\begin{array}{l}
e_{\tau}(2n+1)=\sum\limits_{d=0}^r\sum\limits_{j=0}^{2n}(e_{\pi^d}(j)-e_{\pi^{d-1}}(j))e_{\sigma^d}(2n-j)+
\sum\limits_{d=0}^r\sum\limits_{j=0}^{2n}(o_{\pi^d}(j)-o_{\pi^{d-1}}(j))o_{\sigma^d}(2n-j),\\
\\
e_{\tau}(2n)=\sum\limits_{d=0\ }^r\sum\limits_{j=0,2,4,\ldots,2n-2}(e_{\pi^d}(j)-e_{\pi^{d-1}}(j))o_{\sigma^d}(2n-j-1)+(o_{\pi^d}(j)-o_{\pi^{d-1}}(j))e_{\sigma^d}(2n-j-1)\\
\qquad\quad+\sum\limits_{d=0\
}^r\sum\limits_{j=1,3,5,\ldots,2n-1}(e_{\pi^d}(j)-e_{\pi^{d-1}}(j))e_{\sigma^d}(2n-j-1)+(o_{\pi^d}(j)-o_{\pi^{d-1}}(j))o_{\sigma^d}(2n-j-1).
\end{array}$$
\end{proposition}
\begin{proof}
We use induction. Clearly, the result holds for $n=1$. Now let
$\pi\in\SS_n(132)$ such that $\pi_{j+1}=n$, $0\leq j\leq n-1$.
Then $\beta=(\pi_1,\ldots,\pi_j)$ is a $132$-avoiding permutation
on the letters $n-1,n-2,\dots,n-j$ and
$\gamma=(\pi_{j+1},\ldots,\pi_n)$ is a $132$-avoiding permutation
on the letters $n-j-1,n-j-2,\dots,1$. If we assume that $\beta$
avoids $\pi^d$ and contains $\pi^{d-1}$, then $\gamma$ avoids
$\sigma^d$, where $d=0,1,2\dots,r$. Besides,
$$\sign(\pi)=(-1)^{(j+1)(n-j-1)}\sign(\beta)\sign(\gamma)=(-1)^{(j+1)(n-1)}\sign(\beta)\sign(\gamma),$$
equivalently,
$$\sign(\pi)=\left\{\begin{array}{ll}
                       \sign(\beta)\cdot\sign(\gamma),&\mbox{if}\, n\,\mbox{odd}\\
                        (-1)^{j+1}\cdot\sign(\beta)\cdot\sign(\gamma),\quad&\mbox{if}\, n\,\mbox{even}
                        \end{array}\right.$$
Hence, if summing over all $d=0,1,\ldots,r$ and $j=0,1,2\ldots,n$
together with use of the fact that the number of even \resp{odd}
permutations in $\EE_n$ \resp{$\OO_n$} avoiding $\beta$ and
containing $\gamma$ is given by $e_{\beta}(n)-e_{\beta,\gamma}(n)$
\resp{ $o_{\beta}(n)-o_{\beta,\gamma}(n)$} we get the desired
result.
\end{proof}
Similarly to Proposition~\ref{gen1} we obtain the following
result.

\begin{proposition}\label{gen2}
Let $k\geq 1$ and $n\geq 1$, then
$$\begin{array}{l}
o_{\tau}(2n+1)=\sum\limits_{d=0}^r\sum\limits_{j=0}^{2n}(e_{\pi^d}(j)-e_{\pi^{d-1}}(j))o_{\sigma^d}(2n-j)+
\sum\limits_{d=0}^r\sum\limits_{j=0}^{2n}(o_{\pi^d}(j)-o_{\pi^{d-1}}(j))e_{\sigma^d}(2n-j),\\
\\
o_{\tau}(2n)=\sum\limits_{d=0\ }^r\sum\limits_{j=0,2,4,\ldots,2n-2}(e_{\pi^d}(j)-e_{\pi^{d-1}}(j))e_{\sigma^d}(2n-j-1)+(o_{\pi^d}(j)-o_{\pi^{d-1}}(j))o_{\sigma^d}(2n-j-1)\\
\qquad\quad+\sum\limits_{d=0\
}^r\sum\limits_{j=1,3,5,\ldots,2n-1}(e_{\pi^d}(j)-e_{\pi^{d-1}}(j))o_{\sigma^d}(2n-j-1)+(o_{\pi^d}(j)-o_{\pi^{d-1}}(j))e_{\sigma^d}(2n-j-1).
\end{array}$$
\end{proposition}

Our present aim is to find the generating functions $\EE_\tau(x)$
and $\OO_\tau(x)$; thus we need the following lemma which holds
immediately by definitions.

\begin{lemma}\label{fact1}
let $\{a_n\}_{n\geq 0}$ and $\{b_n\}_{n\geq 0}$ be two sequences,
and the corresponding generating functions are $a(x)$ and $b(x)$;
respectively. Then
\begin{enumerate}
\item $\sum\limits_{n\geq 0} a_{2n}x^{2n} =\frac{1}{2}
(a(x)+a(-x))$;

\item $\sum\limits_{n\geq 1} a_{2n-1}x^{2n-1}=\frac{1}{2}
(a(x)-a(-x))$;

\item $\sum\limits_{n\geq 1\ } \sum\limits_{j=0,2,4,\dots,2n-2}
a_jb_{2n-1-j}x^{2n-1}=\frac{1}{4}(a(x)+a(-x))(b(x)-b(-x))$;

\item $\sum\limits_{n\geq 1\ } \sum\limits_{j=1,3,5,\dots,2n-1} a_jb_{2n-1-j}x^{2n-1}=\frac{1}{4}(a(x)-a(-x))(b(x)+b(-x))$.
\end{enumerate}
\end{lemma}

\begin{theorem}\label{thgen} For any nonempty pattern $\tau\in \SS_k(132)$,
the generating functions $\EE_\tau(x)$ and $\OO_\tau(x)$ are a
rational functions in $x$ satisfying the relations
\begin{equation}
\begin{array}{l}
\EE_{\tau}(x)-\EE_{\tau}(-x)=x\sum\limits_{d=0}^r (\EE_{\pi^d}(x)-\EE_{\pi^{d-1}}(x))\EE_{\sigma^d}(x)+(\EE_{\pi^d}(-x)-\EE_{\pi^{d-1}}(-x))\EE_{\sigma^d}(-x) \\
\qquad\qquad\qquad\quad\,+x\sum\limits_{d=0}^r(\OO_{\pi^d}(x)-\OO_{\pi^{d-1}}(x))\OO_{\sigma^d}(x)+(\OO_{\pi^d}(-x)-\OO_{\pi^{d-1}}(-x))\OO_{\sigma^d}(-x)),
\end{array}\label{eq1}
\end{equation}

\begin{equation}
\begin{array}{l}
\OO_{\tau}(x)-\OO_{\tau}(-x)=x\sum\limits_{d=0}^r (\EE_{\pi^d}(x)-\EE_{\pi^{d-1}}(x))\OO_{\sigma^d}(x)+(\EE_{\pi^d}(-x)-\EE_{\pi^{d-1}}(-x))\OO_{\sigma^d}(-x) \\
\qquad\qquad\qquad\quad\,+x\sum\limits_{d=0}^r(\OO_{\pi^d}(x)-\OO_{\pi^{d-1}}(x))\EE_{\sigma^d}(x)+(\OO_{\pi^d}(-x)-\OO_{\pi^{d-1}}(-x))\EE_{\sigma^d}(-x)),
\end{array}\label{fq1}
\end{equation}

\begin{equation}
\begin{array}{l}
\EE_\tau(x)+\EE_\tau(-x)-2=\\
\quad\quad\qquad=\dfrac{x}{2}\sum\limits_{d=0}^r (\EE_{\pi^d}(x)+\EE_{\pi^d}(-x)-\EE_{\pi^{d-1}}(x)-\EE_{\pi^{d-1}}(-x))(\OO_{\sigma^d}(x)-\OO_{\sigma^d}(-x))\\
\quad\quad\qquad+\dfrac{x}{2}\sum\limits_{d=0}^r (\OO_{\pi^d}(x)+\OO_{\pi^d}(-x)-\OO_{\pi^{d-1}}(x)-\OO_{\pi^{d-1}}(-x))(\EE_{\sigma^d}(x)-\EE_{\sigma^d}(-x))\\
\quad\quad\qquad+\dfrac{x}{2}\sum\limits_{d=0}^r (\EE_{\pi^d}(x)-\EE_{\pi^d}(-x)-\EE_{\pi^{d-1}}(x)+\EE_{\pi^{d-1}}(-x))(\EE_{\sigma^d}(x)+\EE_{\sigma^d}(-x))\\
\quad\quad\qquad+\dfrac{x}{2}\sum\limits_{d=0}^r(\OO_{\pi^d}(x)-\OO_{\pi^d}(-x)-\OO_{\pi^{d-1}}(x)+\OO_{\pi^{d-1}}(-x))(\OO_{\sigma^d}(x)+\OO_{\sigma^d}(-x)),
\end{array}\label{eq2}
\end{equation}
and
\begin{equation}
\begin{array}{l}
\OO_\tau(x)+\OO_\tau(-x)=\\
\quad\quad\qquad=\dfrac{x}{2}\sum\limits_{d=0}^r (\EE_{\pi^d}(x)+\EE_{\pi^d}(-x)-\EE_{\pi^{d-1}}(x)-\EE_{\pi^{d-1}}(-x))(\EE_{\sigma^d}(x)-\EE_{\sigma^d}(-x))\\
\quad\quad\qquad+\dfrac{x}{2}\sum\limits_{d=0}^r (\OO_{\pi^d}(x)+\OO_{\pi^d}(-x)-\OO_{\pi^{d-1}}(x)-\OO_{\pi^{d-1}}(-x))(\OO_{\sigma^d}(x)-\OO_{\sigma^d}(-x))\\
\quad\quad\qquad+\dfrac{x}{2}\sum\limits_{d=0}^r (\EE_{\pi^d}(x)-\EE_{\pi^d}(-x)-\EE_{\pi^{d-1}}(x)+\EE_{\pi^{d-1}}(-x))(\OO_{\sigma^d}(x)+\OO_{\sigma^d}(-x))\\
\quad\quad\qquad+\dfrac{x}{2}\sum\limits_{d=0}^r(\OO_{\pi^d}(x)-\OO_{\pi^d}(-x)-\OO_{\pi^{d-1}}(x)+\OO_{\pi^{d-1}}(-x))(\EE_{\sigma^d}(x)+\EE_{\sigma^d}(-x)).
\end{array}\label{fq2}
\end{equation}
\end{theorem}
\begin{proof}
Using Propositions~\ref{gen1} and \ref{gen2} together with
Lemma~\ref{fact1} we get Equations~\ref{eq1}-\ref{fq2}.
Rationalities of $\EE_\tau(x)$ and $\OO_\tau(x)$ for $\tau\neq\vr$
follows easily by induction.
\end{proof}

As a remark, the above theorem holds for the empty pattern without
rationalities of $\EE_\tau(x)$ and $\OO_\tau(x)$.

Using Theorem~\ref{thgen} we get the main result of~\cite{MV2}.

\begin{corollary}\label{ffc}{\rm(see Mansour and
Vainshtein~\cite[Theorem~2.1]{MV2})} For any nonempty pattern
$\tau\in \SS_k(132)$, the generating function $\F_\tau(x)$ is a
rational function in $x$ satisfying the relation
$$\F_{\tau}(x)=1+x\sum\limits_{d=0}^r(\F_{\pi^d}(x)-\F_{\pi^{d-1}}(x))\F_{\sigma^d}(x).$$
\end{corollary}
\begin{proof}
If adding Equations~\ref{eq1} and \ref{fq1} we have
\begin{equation}
\F_{\tau}(x)-\F_{\tau}(-x)=x\sum\limits_{d=0}^r(\F_{\pi^d}(x)-\F_{\pi^{d-1}}(x))\F_{\sigma^d}(x)+(\F_{\pi^d}(-x)-\F_{\pi^{d-1}}(-x))\F_{\sigma^d}(-x).
\label{meqa}\end{equation} If adding Equations~\ref{eq2} and
\ref{fq2} we get
$$\begin{array}{l}
\F_\tau(x)+\F_\tau(-x)-2=\\
\quad\quad\qquad=\dfrac{x}{2}\sum\limits_{d=0}^r (\F_{\pi^d}(x)+\F_{\pi^d}(-x)-\F_{\pi^{d-1}}(x)-\F_{\pi^{d-1}}(-x))(\F_{\sigma^d}(x)-\F_{\sigma^d}(-x))\\
\quad\quad\qquad+\dfrac{x}{2}\sum\limits_{d=0}^r(\F_{\pi^d}(x)-\F_{\pi^d}(-x)-\F_{\pi^{d-1}}(x)+\F_{\pi^{d-1}}(-x))(\F_{\sigma^d}(x)+\F_{\sigma^d}(-x)).
\end{array}$$
equivalently,
\begin{equation}
\F_{\tau}(x)+\F_{\tau}(-x)-2=x\sum\limits_{d=0}^r(\F_{\pi^d}(x)-\F_{\pi^{d-1}}(x))\F_{\sigma^d}(x)-(\F_{\pi^d}(-x)-\F_{\pi^{d-1}}(-x))\F_{\sigma^d}(-x).
\label{meqb}\end{equation} Hence, if adding Equations~\ref{meqa}
and \ref{meqb} then we obtain that
$$\F_{\tau}(x)=1+x\sum\limits_{d=0}^r(\F_{\pi^d}(x)-\F_{\pi^{d-1}}(x))\F_{\sigma^d}(x).$$
Rationality of $\F_\tau(x)$ for $\tau\neq\vr$ follows easily by
rationalities of $\OO_\tau(x)$ and $\EE_\tau(x)$ (see
Theorem~\ref{thgen} and Equation~\ref{eq3}).
\end{proof}

Our present aim is to find explicitly the generating functions
$\EE_\tau(x)$ and $\OO_\tau(x)$ for several cases of $\tau$; thus
we need the following notation. We denote the generating function
$\EE_\tau(x)-\OO_\tau(x)$ by $\MM_\tau(x)$ for any pattern $\tau$.

\begin{theorem}\label{genmm}
For any $\tau\in\SS_k(132)$,
\begin{equation}
\MM_{\tau}(x)-\MM_{\tau}(-x)=x\sum\limits_{d=0}^r(\MM_{\pi^d}(x)-\MM_{\pi^{d-1}}(x))\MM_{\sigma^d}(x)+(\MM_{\pi^d}(-x)-\MM_{\pi^{d-1}}(-x))\MM_{\sigma^d}(-x),
\label{mmq1}
\end{equation}
and
\begin{equation}
\MM_\tau(x)+\MM_\tau(-x)-2=x\sum\limits_{d=0}^r(\MM_{\pi^d}(x)-\MM_{\pi^{d-1}}(x))\MM_{\sigma^d}(-x)-(\MM_{\pi^d}(-x)-\MM_{\pi^{d-1}}(-x))\MM_{\sigma^d}(x).\label{mmq2}
\end{equation}
\end{theorem}
\begin{proof}
If subtracting Equation~\ref{fq1} from Equation~\ref{eq1} then we
get Equation~\ref{mmq1}, and if subtracting Equation~\ref{fq2}
from Equation~\ref{eq2} then we get Equation~\ref{mmq2}.
\end{proof}

\begin{corollary}\label{mmc}
Let $\tau=(\beta,k)\in\SS_k(132)$. Then
$$\MM_\tau(x)=\frac{2(1+x\MM_\beta(-x))}{(1-x\MM_\beta(x))^2+(1+x\MM_\beta(-x))^2}.$$
\end{corollary}
\begin{proof}
Equation~\ref{mmq1} for $\tau=(\beta,k)$ yields
    $$(1-x\MM_\beta(x))\MM_\tau(x)-(1+x\MM_\beta(-x))\MM_\tau(-x)=0,$$
and Equation~\ref{mmq2} for $\tau=(\beta,k)$ yields
    $$(1+x\MM_\beta(-x))\MM_\tau(x)+(1-x\MM_\beta(x))\MM_\tau(-x)=2.$$
Hence, the rest is easy to check by the above two equations.
\end{proof}

\begin{example}\label{ex1}
Let $\tau=12$ and $\beta=1$. Since $\EE_\beta(x)=1$ and
$\OO_\beta(x)=0$ we get that $\MM_\beta(x)=1$. Corollary~\ref{mmc}
for $\tau=12$ yields
    $$\MM_{12}(x)=\EE_{12}(x)-\OO_{12}(x)=\frac{2(1+x)}{(1-x)^2+(1+x)^2}=\frac{1+x}{1+x^2}.$$
On the other hand, Corollary~\ref{ffc} together with
Equation~\ref{eq3} we have
    $$\F_{12}(x)=\EE_{12}(x)+\OO_{12}(x)=\frac{1}{1-x}.$$
Hence,
    $$\EE_{12}(x)=\frac{1+x}{1-x^4}\mbox{ and }\OO_{12}(x)=\frac{x^2(1+x)}{1-x^4}.$$
\end{example}
\subsection{Pattern $\tau=\vr$} Let us consider the case
$\tau=\vr$ as the first case which examined by Simion and Schmidt
\cite{SS}.

\begin{theorem}\label{thss} We have
    $$\EE(x)=\frac{1}{2}(C(x)+1)+\frac{x}{2}C(x^2)\mbox{ and }\OO(x)=\frac{1}{2}(C(x)-1)-\frac{x}{2}C(x^2).$$
In other words, for all $n\geq1$,
\begin{enumerate}
\item $|\EE_{2n-2}(132)|=\frac{1}{2}C_{2n-2}$;
\item $|\OO_{2n-2}(132)|=\frac{1}{2}C_{2n-2}$;
\item $|\EE_{2n-1}(132)|=\frac{1}{2} (C_{2n-1}+C_{n-1})$;
\item $|\OO_{2n-1}(132)|=\frac{1}{2} (C_{2n-1}-C_{n-1})$.
\end{enumerate}
\end{theorem}
\begin{proof}
First of all, let us define $\MM(x)=\MM_\vr(x)$,
$\F(x)=\F_\vr(x)$, $\EE(x)=\EE_\vr(x)$, and $\OO(x)=\OO_\vr(x)$.
Using the same arguments in the proof of Corollary~\ref{mmc} we
get
$$\left\{\begin{array}{l}
(1-x\MM(x))\MM(x)-(1+x\MM(-x))\MM(-x)=0,\\
(1+x\MM(-x))\MM(x)+(1-x\MM(x))\MM(-x)=2,
\end{array}\right.$$
therefore,
        $$\MM(x)=\EE(x)-\OO(x)=1+xC(x^2).$$
On the other hand, Corollary~\ref{ffc} for $\tau=\vr$ yields
$\F(x)=1+x\F(x)^2$, so by Equation~\ref{eq3} we have that
$$\F(x)=\EE(x)+\OO(x)=C(x).$$
The rest is easy to check.
\end{proof}

\subsection{Pattern $\tau=12\ldots k$}\label{ans1} Let us start be the
following example.
\begin{example}\label{exaa1} {\rm (see Simion and Schmidt~\cite[Proposition~7]{SS})}
Corollary~\ref{mmc} together with Example~\ref{ex1} yield
    $$\MM_{123}(x)=1+x.$$
On the other hand, using the fact that
$\F_{123}(x)=\EE_{123}(x)+\OO_{123}(x)=\frac{1-x}{1-2x}$ {\rm(}see
\cite{SS}{\rm)} we get
$$\EE_{123}(x)=1+x+\frac{x^2}{1-2x}\mbox{ and }\OO_{123}(x)=\frac{x^2}{1-2x}.$$
\end{example}

The case of varying $k$ is more interesting. As an extension of
Example~\ref{exaa1}, let us consider the case $\tau=[k]$, where we
define $[k]=12\ldots k$.

\begin{theorem}\label{thpk}
For all $k\geq 1$,

\noindent{\rm(i)}\,\,
$\EE_{[2k-1]}(x)=\dfrac{1}{2}(R_{2k-1}(x)+xR_{k-1}(x^2)+1)\mbox{
and
}\EE_{[2k]}(x)=\dfrac{1}{2}\left(R_{2k}(x)+\dfrac{(1+xR_k(x^2))R_k(x^2)}{1+x^2R_k^2(x^2)}\right)$,

\noindent{\rm(ii)}
$\OO_{[2k-1]}(x)=\dfrac{1}{2}(R_{2k-1}(x)-xR_{k-1}(x^2)-1)\mbox{
and
}\OO_{[2k]}(x)=\dfrac{1}{2}\left(R_{2k}(x)-\dfrac{(1+xR_k(x^2))R_k(x^2)}{1+x^2R_k^2(x^2)}\right)$.
\end{theorem}
\begin{proof}
We use induction on $k$. Using Example~\ref{ex1} we get that
$\MM_{1}(x)=1$ and $\MM_{12}(x)=\frac{1+x}{1+x^2}$. Now, let us
fix $k$ and assume that
$$\MM_{[2k-1]}(x)=1+xR_{k-1}(x^2)\mbox{ and }
\MM_{[2k]}(x)=\frac{(1+xR_k(x^2))R_k(x^2)}{1+x^2R_k^2(x^2)}.$$
Therefore, by the hypothesis of the induction and
Corollary~\ref{mmc} we get

$\begin{array}{l}
\MM_{[2k+1]}=\dfrac{2\left(1+\frac{x(1-xR_k(x^2))R_k(x^2)}{1+x^2R_k^2(x^2)}\right)}
{\left(1-\frac{x(1+xR_k(x^2))R_k(x^2)}{1+x^2R_k^2(x^2)}\right)^2
+\left(1+\frac{x(1-xR_k(x^2))R_k(x^2)}{1+x^2R_k^2(x^2)}\right)^2}\\
\\
\qquad\quad\,\,\,\,=\dfrac{\frac{2(1+xR_k(x^2))}{1+x^2R_k^2(x^2)}}
{\frac{(1-xR_k(x^2))^2}{(1+x^2R_k^2(x^2))^2}+\frac{(1+xR_k(x^2))^2}{(1+x^2R_k^2(x^2))^2}}
=1+xR_k(x^2).
\end{array}$

Also, by the induction hypothesis, Corollary~\ref{mmc}, and
Identity~\ref{rk} we have that

$\begin{array}{l} \MM_{[2k+2]}=\dfrac{2(1+x(1-xR_k(x^2)))}
{\left(1-x(1+xR_k(x^2))\right)^2+\left(1+x(1-xR_k(x^2))\right)^2}\\
\\
\qquad\quad\,\,\,\,=\dfrac{2\left(x+\frac{1}{R_{k+1}(x^2)}\right)}
{\left(\frac{1}{R_{k+1}(x^2)}-x\right)^2+\left(\frac{1}{R_{k+1}(x^2)}+x\right)^2}\\
\\
\qquad\quad\,\,\,\,=\dfrac{2(1+xR_{k+1}(x^2))R_{k+1}(x^2)}
{(1-xR_{k+1}(x^2))^2+(1+xR_{k+1}(x^2))^2}=\dfrac{(1+xR_{k+1}(x^2))R_{k+1}(x^2)}{1+x^2R_{k+1}^2(x^2)}.
\end{array}$

Hence, for all $k\geq 1$,
$$\MM_{[2k-1]}(x)=1+xR_{k-1}(x^2)\mbox{ and }\MM_{[2k]}(x)=\frac{(1+xR_k(x^2))R_k(x^2)}{1+x^2R_k^2(x^2)}.$$
On the other hand, in \cite{CW} (see also \cite{MV1,Kr,MV2,MV4})
was proved
        $$\F_{[m]}(x)=\EE_{[m]}(x)+\OO_{[m]}(x)=R_m(x),$$
for all $m\geq1$. Hence, by the above two equations we get the
desired result.
\end{proof}

\begin{example} Theorem~\ref{thpk} for $k=3,4,5$ yields
$$\begin{array}{l}
    e_{[5]}(n)-o_{[5]}(n)=\frac{1}{2} \left( 1+(-1)^{n+1} \right); \\
    e_{[7]}(n)-o_{[7]}(n)=\sqrt{2}^{n-5} \left( 1+(-1)^{n+1} \right); \\
    e_{[9]}(n)-o_{[9]}(n)=\frac{F_{n-3}}{2} \left( 1+(-1)^{n+1} \right); \\
\end{array}$$
where $F_{n-3}$ is the $(n-3)$th Fibonacci number.
\end{example}

\subsection{Pattern $\tau=213\ldots k$} Let us start by the
following example.
\begin{example}\label{exaa2} {\rm (see Simion and Schmidt~\cite[Proposition~7]{SS})}
Corollary~\ref{mmc} for $\tau=213$ together with the fact that
$\MM_{21}(x)=\EE_{21}(x)=\frac{1}{1-x}$ and $\OO_{21}(x)=0$ yield
$$\EE_{213}(x)=\frac{(1-x)(1-4x^2+4x^4)}{(1-2x)(1-3x^2+4x^4)},\quad
\OO_{213}(x)=\frac{(1-x)x^2}{(1-2x)(1-3x^2+4x^4)}.$$
\end{example}

The case of varying $k$ is more interesting. As an extension of
Example~\ref{exaa2}, by using Corollary~\ref{mmc}, and induction
on $k$ (Similarly to Theorem~\ref{thpk}) we get

\begin{theorem}\label{thpkk}
For all $k\geq1$,
$$\begin{array}{l}
\MM_{2134\ldots(2k)}(x)=\dfrac{(1+2x)\left(U_{2k-1}\sx-U_{2k}\sx+\frac{x^2+2x-1}{1-3x^2}\right)}{U_{2k}\sx-2xU_{2k+1}\sx+\frac{4x^4}{1-3x^2}},\\
\\
$$\MM_{2134\ldots(2k-1)}(x)=\dfrac{(1+2x)\left(U_{2k-2}\sx-U_{2k-1}\sx+\frac{x^2+2x-1}{1-3x^2}\right)}{U_{2k-1}\sx-2xU_{2k}\sx-\frac{2x(1-5x^2)}{1-3x^2}}.
\end{array}$$
\end{theorem}

By Theorem~\ref{thpkk} together with use of the fact that
$\F_{2134\ldots k}(x)=R_k(x)$ (see \cite[Theorem~2.6]{MV2}) we get

\begin{corollary}\label{ccpkk}
For all $k\geq1$,

{\rm(i)}
$\EE_{2134\ldots(2k)}(x)=\frac{1}{2}\left[\frac{U_{2k-1}\txx}{\sqrt{x}U_{2k}\txx}+\frac{(1+2x)\left(U_{2k-1}\sx-U_{2k}\sx+\frac{x^2+2x-1}{1-3x^2}\right)}{U_{2k}\sx-2xU_{2k+1}\sx+\frac{4x^4}{1-3x^2}}\right]$;

{\rm(ii)}
$\OO_{2134\ldots(2k)}(x)=\frac{1}{2}\left[\frac{U_{2k-1}\txx}{\sqrt{x}U_{2k}\txx}-\frac{(1+2x)\left(U_{2k-1}\sx-U_{2k}\sx+\frac{x^2+2x-1}{1-3x^2}\right)}{U_{2k}\sx-2xU_{2k+1}\sx+\frac{4x^4}{1-3x^2}}\right]$;

{\rm(iii)}
$\EE_{2134\ldots(2k-1)}(x)=\frac{1}{2}\left[\frac{U_{2k-2}\txx}{\sqrt{x}U_{2k-1}\txx}+\frac{(1+2x)\left(U_{2k-2}\sx-U_{2k-1}\sx+\frac{x^2+2x-1}{1-3x^2}\right)}{U_{2k-1}\sx-2xU_{2k}\sx-\frac{2x(1-5x^2)}{1-3x^2}}\right]$;

{\rm(iv)}
$\OO_{2134\ldots(2k-1)}(x)=\frac{1}{2}\left[\frac{U_{2k-2}\txx}{\sqrt{x}U_{2k-1}\txx}-\frac{(1+2x)\left(U_{2k-2}\sx-U_{2k-1}\sx+\frac{x^2+2x-1}{1-3x^2}\right)}{U_{2k-1}\sx-2xU_{2k}\sx-\frac{2x(1-5x^2)}{1-3x^2}}\right]$.
\end{corollary}

\subsection{Pattern $(d+1)(d+2)\ldots k12\ldots d$} In this
subsection we consider the case $\tau=[k,d]$ where
$[k,d]=(d+1)(d+2)\ldots k12\ldots d$. Following to
Theorem~\ref{thpk}, our present aim is to find explicitly the
generating functions $\EE_\tau(x)$ and $\OO_\tau(x)$ where
$\tau=[k,d]$; thus we need to consider four cases either $k$ even
or odd, and either $d$ even or odd. First of all, using
Theorem~\ref{genmm} for $\tau=[k,d]$ we state the following fact.

\begin{lemma}\label{gdd}
Let $k\geq 2$, $1\leq d\leq k-1$, and $\tau=[k,d]$. Then
$$\left\{\begin{array}{l}
(1+x\MM_{[k-d-1]}(-x)-x\MM_{[d]}(-x))\MM_{\tau}(x)+(1-x\MM_{[k-d-1]}(x)+x\MM_{[d]}(x))\MM_{\tau}(-x)\\
\quad\qquad\qquad\qquad\qquad\qquad\qquad\qquad=2+x\MM_{[k-d-1]}(-x)\MM_{[d]}(x)-x\MM_{[k-d-1]}(x)\MM_{[d]}(-x),\\
\\
(1-x\MM_{[k-d-1]}(x)-x\MM_{[d]}(x))\MM_{\tau}(x)-(1+x\MM_{[k-d-1]}(-x)+x\MM_{[d]}(-x))\MM_{\tau}(-x)\\
\quad\qquad\qquad\qquad\qquad\qquad\qquad\qquad=-x\MM_{[k-d-1]}(x)\MM_{[d]}(x)-x\MM_{[k-d-1]}(-x)\MM_{[d]}(-x).
\end{array}\right.$$
\end{lemma}

\subsubsection{$k$ and $d$ are odd numbers}
Now, we ready to consider the first case $k$ and $d$ are odd
numbers.

\begin{theorem}\label{ths1}
Let $0\leq d\leq k$. Then
$$\MM_{[2k+1,2d+1]}(x)=1+xR_k(x^2).$$
\end{theorem}
\begin{proof}
Theorem~\ref{thpk} yields $\MM_{2k-2d-1}(x)=1+xR_{k-d-1}(x^2)$ and
$\MM_{2d+1}(x)=1+xR_d(x^2)$. Therefore, by use of Lemma~\ref{gdd}
for $\tau=[2k+1,2d+1]$ we get
$$\left\{\begin{array}{l}
\MM_{\tau}(x)+\MM_{\tau}(-x)=2,\\
\\
(1-2x-x^2R_d(x^2)-x^2R_{k-d-1}(x^2))\MM_{\tau}(x)-(1+2x-x^2R_{d}(x^2)-x^2R_{k-d-1}(x^2))\MM_{\tau}(-x)\\
\qquad\qquad\qquad\qquad\qquad\qquad\qquad\qquad\qquad\qquad\qquad\qquad\qquad\qquad=-2x-2x^3R_d(x^2)R_{k-d-1}(x^2),
\end{array}\right.$$
so,
$$\MM_\tau(x)=1+\frac{x(1-x^2R_d(x^2)R_{k-d-1}(x^2))}{1-x^2(R_d(x^2)+R_{k-d-1}(x^2))}.$$
By using the following identities (see~\cite{MV2})
\begin{equation}
1-x^2R_p(x^2)R_q(x^2)=\frac{U_{p+q}\sx}{U_p\sx U_q\sx}\mbox{ and }
1-x^2(R_p(x^2)+R_q(x^2))=\frac{xU_{p+q+1}\sx}{U_p\sx U_q\sx},
\label{irks}
\end{equation}
we have
$$\MM_\tau(x)=1+\frac{U_{k-1}\sx}{U_k\sx},$$
and by Identity~\ref{drk} we get the desired result.
\end{proof}

By Theorem~\ref{ths1} together with use of the equation
$\F_{[k,d]}(x)=R_k(x)$ (see \cite[Theorem~2.4]{MV2}) we get

\begin{corollary}\label{ccc1}
For all $0\leq d\leq k$,
$$\EE_{[2k+1,2d+1]}(x)=\frac{1}{2}(R_{2k+1}(x)+xR_{k}(x^2)+1)\mbox{
and }
\OO_{[2k+1,2d+1]}(x)=\frac{1}{2}(R_{2k+1}(x)-xR_{k}(x^2)-1).$$
\end{corollary}

A comparison of Corollary~\ref{ccc1} for values of $d$ suggests
that there should exist a bijection between the sets
$\EE_n(132,2\ldots(2k+1)1)$ and
$\EE_n(132,(2d+2)(2d+3)\ldots(2k+1)12\ldots(2d+1))$ for any $d$.
However, we failed to produce such a bijection, and finding it
remains a challenging open question.

\subsubsection{$k$ odd number and $d$ even number}
Now, let us consider the case where $k$ odd number and $d$ even
number.

\begin{theorem}\label{ths2}
Let $1\leq d\leq k$. Then
        $$\MM_{[2k+1,2d]}(x)=1+xR_k(x^2).$$
\end{theorem}
\begin{proof}
Let $m=k-d$; solving the system equations in Lemma~\ref{gdd} for
$\tau=[2k+1,2d]$ together with use of Theorem~\ref{thpk} we get
$$\MM_{[2k+1,2d]}(x)=1+\frac{x(R_m(x^2)+R_d(x^2)-R_m(x^2)R_d(x^2))}{1-x^2R_m(x^2)R_d(x^2)}.$$
Using Identities~\ref{irks} we get the desired result.
\end{proof}

By Theorem~\ref{ths2} together with use of the equation
$\F_{[k,d]}(x)=R_k(x)$ (see \cite[Theorem~2.4]{MV2}) we get

\begin{corollary}\label{ccc2}
For all $1\leq d\leq k$,
$$\EE_{[2k+1,2d]}(x)=\dfrac{1}{2}(R_{2k+1}+1+xR_k(x^2))\mbox{ and }
\OO_{[2k+1,2d]}(x)=\dfrac{1}{2}(R_{2k+1}-1-xR_k(x^2).$$
\end{corollary}

A comparison of Corollary~\ref{ccc1} with Corollary~\ref{ccc2}
suggests that there should exist a bijection between the sets
$\EE_n(132,12\ldots(2k+1))$ and $\EE_n(132,[2k+1,d])$ for any $d$.
However, we failed to produce such a bijection, and finding it
remains a challenging open question.

\subsubsection{$k$ even number and $d$ odd number}
Similarly as above subsections, we can consider the case where $k$
even number and $d$ odd number.

\begin{theorem}\label{ths3}
Let $0\leq d\leq k-1$ and $m=k-d-1$. Then the generating function
$\MM_{[2k,2d+1]}(x)$ is given by
$$\frac{\bigl(1-x^2(R_m(x^2)+R_d(x^2))+x(1-x^2R_m(x^2)R_d(x^2))\bigr)(1+x^2R_m(x^2)R_d(x^2))}
{1-x^2(1+R_m^2(x^2))(1+x^2R_m^2(x^2))}.$$
\end{theorem}

By Theorem~\ref{ths3} together with use of the equation
$\F_{[k,d]}(x)=R_k(x)$ (see \cite[Theorem~2.4]{MV2}) we get

\begin{corollary}\label{ccc3}
Let $0\leq d\leq k-1$ and $m=k-d-1$.

{\rm(i)} The generating function $\EE_{[2k,2d+1]}(x)$ is given by
$$\dfrac{1}{2}\left(R_{2k}(x)+\frac{\bigl(1-x^2(R_m(x^2)+R_d(x^2))+x(1-x^2R_m(x^2)R_d(x^2))\bigr)(1+x^2R_m(x^2)R_d(x^2))}
{1-x^2(1+R_m^2(x^2))(1+x^2R_m^2(x^2))}\right).$$ {\rm(ii)} The
generating function $\OO_{[2k,2d+1]}(x)$ is given by
$$\dfrac{1}{2}\left(R_{2k}(x)-\frac{\bigl(1-x^2(R_m(x^2)+R_d(x^2))+x(1-x^2R_m(x^2)R_d(x^2))\bigr)(1+x^2R_m(x^2)R_d(x^2))}
{1-x^2(1+R_m^2(x^2))(1+x^2R_m^2(x^2))}\right).$$
\end{corollary}

\subsubsection{$k$ and $d$ even numbers}
Similarly as above subsections, we can consider the case where $k$
and $d$ even numbers.

\begin{theorem}\label{ths4}
Let $1\leq d\leq k-1$ and $m=k-d-1$. Then the generating function
$\MM_{[2k,2d]}(x)$ is given by
$$\frac{1}{x}-\frac{\bigl(1-x^2(R_m(x^2)-R_d(x^2))\bigr)\bigl(1-x^2(R_d(x^2)+R_m(x^2))-x(1-x^2R_d(x^2)R_m(x^2))\bigr)}
{x+x^3(1+x^2R_d^2(x^2))(1-2R_m(x^2)+x^2R_m^2(x^2))}.$$
\end{theorem}

By Theorem~\ref{ths4} together with use of the equation
$\F_{[k,d]}(x)=R_k(x)$ (see \cite[Theorem~2.4]{MV2}) we get

\begin{corollary}\label{ccc4}
Let $1\leq d\leq k-1$ and $m=k-d-1$.

{\rm(i)} The generating function $\EE_{[2k,2d+1]}(x)$ is given by
$$\frac{1}{2}\biggl(R_{2k}(x)+\frac{1}{x}-\frac{\bigl(1-x^2(R_m(x^2)-R_d(x^2))\bigr)\bigl(1-x^2(R_d(x^2)+R_m(x^2))-x(1-x^2R_d(x^2)R_m(x^2))\bigr)}
{x+x^3(1+x^2R_d^2(x^2))(1-2R_m(x^2)+x^2R_m^2(x^2))}.$$

{\rm(ii)} The generating function $\OO_{[2k,2d+1]}(x)$ is given by
$$\frac{1}{2}\biggl(R_{2k}(x)-\frac{1}{x}+\frac{\bigl(1-x^2(R_m(x^2)-R_d(x^2))\bigr)\bigl(1-x^2(R_d(x^2)+R_m(x^2))-x(1-x^2R_d(x^2)R_m(x^2))\bigr)}
{x+x^3(1+x^2R_d^2(x^2))(1-2R_m(x^2)+x^2R_m^2(x^2))}\biggr).$$
\end{corollary}

\subsection{Wedge patterns}
For a further generalization of the results in the pervious
subsections, consider the following definition. We say that
$\tau\in\SS_k$ is a {\it wedge\/} pattern if it can be represented
as $\tau=(\tau^1,\rho^1,\dots,\tau^r,\rho^r)$ so that each of
$\tau^i$ is nonempty, $(\rho^1,\rho^2,\dots,\rho^r)$ is a layered
permutation of $1,\dots,s$ for some $s$, and
$(\tau^1,\tau^2,\dots,\tau^r)=(s+1,s+2,\dots,k)$. For example,
$645783912$ and $456378129$ are wedge patterns. Evidently, $[k,d]$
is a wedge pattern for any $d$. We say that $\tau\in\SS_k$ is an
{\it odd-wedge\/} pattern if it a wedge pattern such that the
length of $(\tau^1,\rho^1,\ldots,\tau^p,\rho^p)$ is given by odd
number for all $p=1,2,\ldots,r$. For example, $23145$ and $34251$
are odd-wedge patterns. Evidently, $[2k+1,d]$ is an odd-wedge
pattern for any $d$.

\begin{theorem}\label{wedge}
$\MM_\tau(x)=1+xR_k(x^2)$ for any odd-wedge pattern
$\tau\in\SS_{2k+1}(132)$.
\end{theorem}
\begin{proof}
We proceed by induction on $r$. If $r=1$ then $\tau=[2k+1,d]$ for
some $d$, and the result is true by Theorems~\ref{ths1} and
\ref{ths2}. For an arbitrary $r>1$, $\tau$ looks like either
$$\tau=(\tau',2p+2d+2,2p+2d+3,\ldots,2k+1,1,2,\ldots,2d),$$
or
$$\tau=(\tau',2p+2d+3,2p+2d+4,\ldots,2k+1,1,2,\ldots,2d+1),$$
for some $d$ and $p$, where
$\tau'=(\tau^1,\rho^1,\ldots,\tau^{r-1},\rho^{r-1})$.

The first case; $\tau'$ contains $2p+1$ elements and it is an
odd-wedge pattern, so by induction $\MM_{\tau'}(x)=1+xR_p(x^2)$,
so Corollary~\ref{mmc} gives
$$\MM_{(\tau',2p+2d+2)}(x)=\frac{(1+xR_{p+1}(x^2))R_{p+1}(x^2)}{1+x^2R_{p+1}^2(x^2)},$$
and then
$$\MM_{(\tau',2p+2d+2,2p+2d+3)}(x)=1+xR_{p+1}(x^2).$$
Therefore, by induction we have
$$\MM_{(\tau',2p+2d+2,2p+2d+3,\ldots,2k)}(x)=\frac{(1+xR_{k-d}(x^2))R_{k-d}(x^2)}{1+x^2R_{k-d}^2(x^2)}.$$
Hence, Theorem~\ref{genmm} for
$\tau=(\tau',2p+2d+2,2p+2d+3,\ldots,2k+1,1,2,\ldots,2d)$ yields
(similarly as Theorem~\ref{ths2})
$$\MM_\tau(x)=1+\frac{x(R_{k-d}(x^2)+R_d(x^2)-R_{k-d}(x^2)R_d(x^2))}{1-x^2R_{k-d}(x^2)R_d(x^2)},$$
and using Identities~\ref{irks} we have that
$\MM_\tau(x)=1+xR_k(x^2)$.

The second case; similarly as the first case we have that
$$\MM_{(\tau',2p+2d+3,2p+2d+4,\ldots,2k)}(x)=1+xR_{k-d-1}(x^2).$$
Hence, Theorem~\ref{genmm} for
$\tau=(\tau',2p+2d+3,2p+2d+4,\ldots,2k+1,1,2,\ldots,2d+1)$ yields
(similarly as Theorem~\ref{ths1})
$$\MM_\tau(x)=1+\frac{x(1-x^2R_d(x^2)R_{k-d-1}(x^2))}{1-x^2(R_d(x^2)+R_{k-d-1}(x^2))},$$
and using Identities~\ref{irks} we have that
$\MM_\tau(x)=1+xR_k(x^2)$.
\end{proof}

A comparison of Theorem~\ref{thpk} with Theorem~\ref{wedge}
suggests that there should exist a bijection between the sets
$\EE_n(132,12\ldots(2k+1))$ \resp{$\OO_n(132,12\ldots(2k+1))$} and
$\EE_n(132,\tau)$ \resp{$\OO_n(132,\tau)$} for any odd-wedge
pattern $\tau$. However, we failed to produce such a bijection,
and finding it remains a challenging open question.

\begin{corollary}\label{cwedge}
For any odd-wedge pattern $\tau\in\SS_{2k+1}(132)$,
$$\EE_{\tau}(x)=\dfrac{1}{2}(R_{2k+1}(x)+xR_k(x^2)+1)\mbox{
and }\OO_{\tau}(x)=\dfrac{1}{2}(R_{2k+1}(x)-xR_k(x^2)-1).$$
\end{corollary}
\section{Containing a pattern exactly once}\label{sec3}
Let $e_{\tau;r}(n)$ \resp{$o_{\tau;r}(n)$} denote the number of
even \resp{odd} permutations in $\SS_n(132)$ that contain
$\tau\in\SS_k(132)$ exactly $r$ times, and $e_{\tau;1}^\rho(n)$
\resp{$o_{\tau;r}^\rho(n)$} denote the number of even \resp{odd}
permutations in $\EE_n(132,\rho)$ that contain $\tau\in\SS_k(132)$
exactly $r$ times. We denote by $\EE_{\tau;r}(x)$ and
$\EE_{\tau;r}^\rho(x)$ \resp{$\OO_{\tau;r}(x)$ and
$\OO_{\tau;r}^\rho(x)$} the corresponding ordinary generating
functions. Using the argument proof of Theorem~\ref{thgen} we get
as follows.

\begin{theorem}\label{thgencc} For any $\tau=(\tau^0,m_0,\dots,\tau^r,m_r)$  be the canonical
decomposition of nonempty $\tau\in S_k(132)$, then
\begin{equation}
\begin{array}{l}
\EE_{\tau;1}(x)-\EE_{\tau;1}(-x)=x\sum\limits_{d=0}^{r+1}\EE_{\pi^{d-1};1}^{\pi^d}(x)\EE_{\sigma^d;1}^{\sigma^{d-1}}(x)+\EE_{\pi^{d-1};1}^{\pi^d}(-x)\EE_{\sigma^d;1}^{\sigma^{d-1}}(-x)+ \\
\qquad\qquad\qquad\qquad\,\,+x\sum\limits_{d=0}^{r+1}\OO_{\pi^{d-1};1}^{\pi^d}(x)\OO_{\sigma^d;1}^{\sigma^{d-1}}(x)+\OO_{\pi^{d-1};1}^{\pi^d}(-x)\OO_{\sigma^d;1}^{\sigma^{d-1}}(-x),
\end{array}\label{cceq1}
\end{equation}

\begin{equation}
\begin{array}{l}
\OO_{\tau;1}(x)-\OO_{\tau;1}(-x)=x\sum\limits_{d=0}^{r+1}\EE_{\pi^{d-1};1}^{\pi^d}(x)\OO_{\sigma^d;1}^{\sigma^{d-1}}(x)+\EE_{\pi^{d-1};1}^{\pi^d}(-x)\OO_{\sigma^d;1}^{\sigma^{d-1}}(-x)+ \\
\qquad\qquad\qquad\qquad\,\,\,\,\,\,+x\sum\limits_{d=0}^{r+1}\OO_{\pi^{d-1};1}^{\pi^d}(x)\EE_{\sigma^d;1}^{\sigma^{d-1}}(x)+\OO_{\pi^{d-1};1}^{\pi^d}(-x)\EE_{\sigma^d;1}^{\sigma^{d-1}}(-x),
\end{array}\label{cceq2}
\end{equation}

\begin{equation}
\begin{array}{l}
\EE_{\tau;1}(x)+\EE_{\tau;1}(-x)=\dfrac{x}{2}\sum\limits_{d=0}^{r+1}(\EE_{\pi^{d-1};1}^{\pi^d}(x)+\EE_{\pi^{d-1};1}^{\pi^d}(-x))(\OO_{\sigma^d;1}^{\sigma^{d-1}}(x)-\OO_{\sigma^d;1}^{\sigma^{d-1}}(-x))+\\
\quad\quad\qquad\qquad\qquad\,\,\,\,\,+\dfrac{x}{2}\sum\limits_{d=0}^{r+1}(\OO_{\pi^{d-1};1}^{\pi^d}(x)+\OO_{\pi^{d-1};1}^{\pi^d}(-x))(\EE_{\sigma^d;1}^{\sigma^{d-1}}(x)-\EE_{\sigma^d;1}^{\sigma^{d-1}}(-x))+\\
\quad\quad\qquad\qquad\qquad\,\,\,\,\,+\dfrac{x}{2}\sum\limits_{d=0}^{r+1}(\EE_{\pi^{d-1};1}^{\pi^d}(x)-\EE_{\pi^{d-1};1}^{\pi^d}(-x))(\EE_{\sigma^d;1}^{\sigma^{d-1}}(x)+\EE_{\sigma^d;1}^{\sigma^{d-1}}(-x))+\\
\quad\quad\qquad\qquad\qquad\,\,\,\,\,+\dfrac{x}{2}\sum\limits_{d=0}^{r+1}(\OO_{\pi^{d-1};1}^{\pi^d}(x)-\OO_{\pi^{d-1};1}^{\pi^d}(-x))(\OO_{\sigma^d;1}^{\sigma^{d-1}}(x)+\OO_{\sigma^d;1}^{\sigma^{d-1}}(-x)),
\end{array}\label{cceq3}
\end{equation}
and
\begin{equation}
\begin{array}{l}
\OO_{\tau;1}(x)+\OO_{\tau;1}(-x)=\dfrac{x}{2}\sum\limits_{d=0}^{r+1}(\EE_{\pi^{d-1};1}^{\pi^d}(x)+\EE_{\pi^{d-1};1}^{\pi^d}(-x))(\EE_{\sigma^d;1}^{\sigma^{d-1}}(x)-\EE_{\sigma^d;1}^{\sigma^{d-1}}(-x))+\\
\quad\quad\qquad\qquad\qquad\,\,\,\,\,+\dfrac{x}{2}\sum\limits_{d=0}^{r+1}(\OO_{\pi^{d-1};1}^{\pi^d}(x)+\OO_{\pi^{d-1};1}^{\pi^d}(-x))(\OO_{\sigma^d;1}^{\sigma^{d-1}}(x)-\OO_{\sigma^d;1}^{\sigma^{d-1}}(-x))+\\
\quad\quad\qquad\qquad\qquad\,\,\,\,\,+\dfrac{x}{2}\sum\limits_{d=0}^{r+1}(\EE_{\pi^{d-1};1}^{\pi^d}(x)-\EE_{\pi^{d-1};1}^{\pi^d}(-x))(\OO_{\sigma^d;1}^{\sigma^{d-1}}(x)+\OO_{\sigma^d;1}^{\sigma^{d-1}}(-x))+\\
\quad\quad\qquad\qquad\qquad\,\,\,\,\,+\dfrac{x}{2}\sum\limits_{d=0}^{r+1}(\OO_{\pi^{d-1};1}^{\pi^d}(x)-\OO_{\pi^{d-1};1}^{\pi^d}(-x))(\EE_{\sigma^d;1}^{\sigma^{d-1}}(x)+\EE_{\sigma^d;1}^{\sigma^{d-1}}(-x)),
\end{array}\label{cceq4}
\end{equation}
\end{theorem}

\begin{remark}
Strictly speaking, Theorem~\ref{thgencc}, unlike
Theorem~\ref{thgen}, is {\it not\/} a recursion for
$\EE_{\tau;1}(x)$ or $\OO_{\tau;1}(x)$, since it involves
functions of type $\EE_{\tau;1}^\rho(x)$ and
$\OO_{\tau;1}^\rho(x)$ (unless $r=0$; see the next subsection).
However, for these functions one can write further recursions
involving similar objects. For example,
$$\begin{array}{l}
\EE_{\pi^{j-1};1}^{pi^j}(x)-\EE_{\pi^{j-1};1}^{\pi^j}(-x)
                            =x\sum\limits_{i=0}^{j}\EE_{\pi^{i-1};1}^{\pi^i}(x)\EE_{\sigma^i_{j-1};1}^{\sigma^{i-1}_{j-1}}(x)+\EE_{\pi^{i-1};1}^{\pi^i}(-x)\EE_{\sigma^i_{j-1};1}^{\sigma^{i-1}_{j-1}}(-x)+ \\
\qquad\qquad\qquad\qquad\qquad\quad\,\,+x\sum\limits_{i=0}^{j}\OO_{\pi^{i-1};1}^{\pi^i}(x)\OO_{\sigma^i_{j-1};1}^{\sigma^{i-1}_{j-1}}(x)+\OO_{\pi^{i-1};1}^{\pi^i}(-x)\OO_{\sigma^i_{j-1};1}^{\sigma^{i-1}_{j-1}}(-x),
\end{array}$$
and
$$\begin{array}{l}
\OO_{\pi^{j-1};1}^{pi^j}(x)-\OO_{\pi^{j-1};1}^{\pi^j}(-x)
                            =x\sum\limits_{i=0}^{j}\EE_{\pi^{i-1};1}^{\pi^i}(x)\OO_{\sigma^i_{j-1};1}^{\sigma^{i-1}_{j-1}}(x)+\EE_{\pi^{i-1};1}^{\pi^i}(-x)\OO_{\sigma^i_{j-1};1}^{\sigma^{i-1}_{j-1}}(-x)+ \\
\qquad\qquad\qquad\qquad\qquad\quad\,\,+x\sum\limits_{i=0}^{j}\OO_{\pi^{i-1};1}^{\pi^i}(x)\EE_{\sigma^i_{j-1};1}^{\sigma^{i-1}_{j-1}}(x)+\OO_{\pi^{i-1};1}^{\pi^i}(-x)\EE_{\sigma^i_{j-1};1}^{\sigma^{i-1}_{j-1}}(-x),
\end{array}$$
where $\sigma_{j-1}^i$ is the $i$th suffix of $\pi^{j-1}$. Though
we have not succeeded to write down a complete set of equations in
the general case (for $\EE_\tau^\rho(x)$ and $\OO_\tau^\rho(x)$),
it is possible to do this in certain particular cases.
\end{remark}

As a corollary of Theorem~\ref{thgencc} we obtain the main result
of \cite{MV2}.

\begin{corollary}\label{ffcc}{\rm(see Mansour and
Vainshtein~\cite[Theorem~3.1]{MV2})} Let
$\tau=(\tau^0,m_0,\dots,\tau^r,m_r)$ be the canonical
decomposition of $\tau\in S_k(132)$, then for all $r\geq0$,
$$G_\tau(x)=x\sum_{j=0}^{r+1}G_{\pi^{j-1}}^{\pi^j}(x)G_{\sigma^j}^{\sigma^{j-1}}(x),$$
where $G_\tau(x)$ is the generating function for the number of
permutations in $\SS_n(132)$ containing $\tau$ exactly once, and
$G_\tau^\rho(x)$ is the generating function for the number of
permutations in $\SS_n(132,\rho)$ containing $\tau$ exactly once.
\end{corollary}
\begin{proof}
If adding Equations~\ref{cceq1} and \ref{cceq2} together with use
of the fact that
$G_\tau^\rho(x)=\EE_\tau^\rho(x)+\OO_\tau^\rho(x)$ for any $\tau$
and $\rho$, we get
\begin{equation}
G_\tau(x)-G_\tau(-x)=x\sum\limits_{d=0}^{r+1}G_{\pi^{d-1}}^{\pi^d}(x)G_{\sigma^d}^{\sigma^{d-1}}(x)+G_{\pi^{d-1}}^{\pi^d}(-x)G_{\sigma^d}^{\sigma^{d-1}}(-x),
\label{hah1}
\end{equation}
and if adding Equations~\ref{cceq3} and \ref{cceq4} we have
\begin{equation}
G_{\tau}(x)+G_{\tau}(-x)=x\sum\limits_{d=0}^{r+1}G_{\pi^{d-1}}^{\pi^d}(x)G_{\sigma^d}^{\sigma^{d-1}}(x)-G_{\pi^{d-1}}^{\pi^d}(-x)G_{\sigma^d}^{\sigma^{d-1}}(-x).
\label{hah2}
\end{equation}
Hence, by adding the Equations~\ref{hah1} and \ref{hah2} we get
the desired result.
\end{proof}

Our present aim is to find explicitly the generating functions
$\EE_{\tau;1}(x)$ and $\OO_{\tau;1}(x)$ for several patterns
$\tau$, thus we the following notations. We define
$\MM_{\tau;r}(x)=\EE_{\tau;r}(x)-\OO_{\tau;r}(x)$ and
$\MM_{\tau;r}^\rho(x)=\EE_{\tau;r}^\rho(x)-\OO_{\tau;r}^\rho(x)$
for any $\tau$ and $\rho$, and $\MM_{\vr;r}^{\rho}(x)=\MM_\rho(x)$
for any $\rho$.

\begin{theorem}\label{gencc}
For any $\tau=(\tau^0,m_0,\dots,\tau^r,m_r)$  be the canonical
decomposition of nonempty $\tau\in S_k(132)$, then
\begin{equation}
\MM_{\tau;1}(x)-\MM_{\tau;1}(-x)=x\sum\limits_{d=0}^{r+1}\MM_{\pi^{d-1};1}^{\pi^d}(x)\MM_{\sigma^d;1}^{\sigma^{d-1}}(x)+\MM_{\pi^{d-1};1}^{\pi^d}(-x)\MM_{\sigma^d;1}^{\sigma^{d-1}}(-x),
\label{ccq1}
\end{equation}
and
\begin{equation}
\MM_{\tau;1}(x)+\MM_{\tau;1}(-x)=x\sum\limits_{d=0}^{r+1}\MM_{\pi^{d-1};1}^{\pi^d}(x)\MM_{\sigma^d;1}^{\sigma^{d-1}}(-x)-\MM_{\pi^{d-1};1}^{\pi^d}(-x)\MM_{\sigma^d;1}^{\sigma^{d-1}}(x).
\label{ccq2}
\end{equation}
\end{theorem}
\begin{proof}
If subtracting Equation~\ref{cceq2} from Equation~\ref{cceq1} then
we get Equation~\ref{ccq1}, and if subtracting
Equation~\ref{cceq4} from Equation~\ref{cceq3} then we get
Equation~\ref{ccq2}.
\end{proof}

\subsection{Pattern $\tau=[k]$} One can try to obtain results similar to
Theorems~\ref{thpk}, \ref{thpkk}, and~\ref{ths1}-\ref{ths4}, but
expressions involved become extremely cumbersome. So we just
consider a simplest wedge pattern, which is the pattern $[k]$.

\begin{theorem}\label{thccca}
For all $m\geq0$,

{\rm(i)} $\MM_{[2m+1];1}(x)=\dfrac{x}{U_m^2\sx}$,

{\rm(ii)}
$\MM_{[2m+2];1}(x)=\dfrac{x^2R_{m+1}^2(x^2)}{(1+x^2R_{m+1}^2(x^2))^2U_m^2\sx}\biggl(1+2xR_{m+1}(x^2)-x^2R_{m+1}^2(x^2)\biggr)$.
\end{theorem}
\begin{proof}
Let $\tau=[k]$, then $r=0$, and it follows from
Theorem~\ref{gencc} that
\begin{equation}
\begin{array}{l}
\MM_{[k];1}(x)-\MM_{[k];1}(-x)=x\MM_{[k-1]}(x)\MM_{[k];1}(x)+x\MM_{[k-1]}(-x)\MM_{[k];1}(-x)\\
\qquad\quad\qquad\qquad\qquad\qquad+x\MM_{[k-1];1}(x)\MM_{[k]}(x)+x\MM_{[k-1];1}(-x)\MM_{[k]}(-x),\\
\\
\MM_{[k];1}(x)+\MM_{[k];1}(-x)=x\MM_{[k-1]}(x)\MM_{[k];1}(-x)-x\MM_{[k-1]}(-x)\MM_{[k];1}(x)\\
\qquad\quad\qquad\qquad\qquad\qquad+x\MM_{[k-1];1}(x)\MM_{[k]}(-x)-x\MM_{[k-1];1}(-x)\MM_{[k]}(x).
\end{array}\label{eqccc}
\end{equation}
Now, let us consider two cases either $k=2m$ or $k=2m+1$ as
follows. (i) Let $k=2m$, Equation~\ref{eqccc} for $k=2m$ and
Theorem~\ref{thpk} together with Identity~\ref{rk} we get
$$\begin{array}{l}
(1-xR_m(x^2))\MM_{[2m];1}(x)-(1+xR_m(x^2))\MM_{[2m];1}(-x)=\\
\qquad\qquad\qquad=\dfrac{x(1+xR_m(x^2))R_m^2(x^2)}{1+x^2R_m^2(x^2)}\MM_{[2m-1];1}(x)+\dfrac{x(1-xR_m(x^2))R_m^2(x^2)}{1+x^2R_m^2(x^2)}\MM_{[2m-1];1}(-x)\\
\\
(1+xR_m(x^2))\MM_{[2m];1}(x)+(1-xR_m(x^2))\MM_{[2m];1}(-x)=\\
\qquad\qquad\qquad=\dfrac{x(1-xR_m(x^2))R_m^2(x^2)}{1+x^2R_m^2(x^2)}\MM_{[2m-1];1}(x)-\dfrac{x(1+xR_m(x^2))R_m^2(x^2)}{1+x^2R_m^2(x^2)}\MM_{[2m-1];1}(-x).
\end{array}$$
If solving the above system of Equations, then using
Identity~\ref{rk} we have
\begin{equation}
\MM_{[2m];1}(x)=\frac{xR_m^2(x^2)}{(1+x^2R_m^2(x^2))^2}\biggl((1-x^2R_m^2(x^2))\MM_{[2m-1];1}(x)-2xR_m(x^2)\MM_{[2m-1];1}(-x)\biggr).
\label{ha1}
\end{equation}

(ii) Let $k=2m+1$, similarly as first case (i) we get
\begin{equation}
\MM_{[2m+1];1}(x)=x\biggl((1-x^2R_m^2(x^2))\MM_{[2m];1}(x)-2xR_m(x^2)\MM_{[2m];1}(-x)\biggr).
\label{ha2}
\end{equation}
If using Equations~\ref{ha1} and \ref{ha2}, then we have that for $m\geq1$,
$$\MM_{[2m+1];1}(x)=x^2R_m^2(x^2)\MM_{[2m-1];1}(x).$$
Besides, by definitions we have that $\MM_{[1];1}(x)=x$, hence
$\MM_{[2m+1];1}(x)=\frac{x}{U_m^2\sx}$. Using Equation~\ref{ha1}
together with the property $U_p^2(t)=U_p^2(-t)$ for all $p$, we
get the desired result.
\end{proof}

Theorem~\ref{thccca} together with the fact that the generating
function for the number permutations in $\SS_n(132)$ containing
$[k]$ exactly once is given by $\frac{1}{U_k^2\txx}$
(see~\cite{MV1}) we have

\begin{corollary}\label{ccra}
For all $m\geq0$,

{\rm(i)}
$\EE_{[2m+1];1}(x)=\frac{1}{2}\left(\frac{1}{U_{2m+1}^2\txx}+\frac{x}{U_m^2\sx}\right)$,

{\rm(ii)}
$\OO_{[2m+1];1}(x)=\frac{1}{2}\left(\frac{1}{U_{2m+1}^2\txx}-\frac{x}{U_m^2\sx}\right)$,

{\rm(iii)}
$\EE_{[2m+2];1}(x)=\frac{1}{2}\left(\frac{1}{U_{2m+2}^2\txx}+\frac{x^2R_{m+1}^2(x^2)}{(1+x^2R_{m+1}^2(x^2))^2U_m^2\sx}\bigl(1+2xR_{m+1}(x^2)-x^2R_{m+1}^2(x^2)\bigr)\right)$,

{\rm(iv)}
$\OO_{[2m+2];1}(x)=\frac{1}{2}\left(\frac{1}{U_{2m+2}^2\txx}-\frac{x^2R_{m+1}^2(x^2)}{(1+x^2R_{m+1}^2(x^2))^2U_m^2\sx}\bigl(1+2xR_{m+1}(x^2)-x^2R_{m+1}^2(x^2)\bigr)\right)$.
\end{corollary}
\section{Furthermore results}\label{sec4}
In this section, we present several directions to generalize and
to extend the results of the previous sections.

\subsection{Statistics on the set $\EE_n(132)$ and on the set $\OO_n(132)$} The first of these
directions is to consider statistics on the set $\EE_n(132)$ (or
on the set $\OO_n(132)$). First of all, let us define
$$\begin{array}{l}
\F(x_1,x_2,\ldots)=\sum_{n\geq0}\sum_{\pi\in\SS_n(132)}\prod_{j\geq1}x_j^{12\ldots j(\pi)},\\
\EE(x_1,x_2,\ldots)=\sum_{n\geq0}\sum_{\pi\in\EE_n(132)}\prod_{j\geq1}x_j^{12\ldots j(\pi)},\\
\OO(x_1,x_2,\ldots)=\sum_{n\geq0}\sum_{\pi\in\OO_n(132)}\prod_{j\geq1}x_j^{12\ldots
j(\pi)},
\end{array}$$
where $12\ldots j(\pi)$ is the number occurrences of the pattern
$12\ldots j$ in $\pi$. We denote the function
$\EE(x_1,x_2,\ldots)-\OO(x_1,x_2,\ldots)$ by
$\MM(x_1,x_2,\ldots)$. Using the same arguments in the proof
Corollary~\ref{mmc} together with the main result of \cite{BCS} we
get as follows.

\begin{theorem}\label{thggg} We have

$$\MM(x_1,x_2,\ldots)=\frac{2(1+x_1\MM(-x_1x_2,x_2x_3,\ldots))}{(1-x_1\MM(x_1x_2,x_2x_3,\ldots))^2+(1+x_1\MM(-x_1x_2,x_2x_3,\ldots))^2},$$
and
$$\F(x_1,x_2,\ldots)=\EE(x_1,x_2,\ldots)+\OO(x_1,x_2,\ldots)=\frac{1}{1-x_1\F(x_1x_2,x_2x_3,\ldots)}.$$
\end{theorem}

An application for Theorem~\ref{thggg} we get the distribution of
the number {\em right to left maxima} on the set $\EE_n(132)$ or
on the set $\OO_n(132)$. Let $\pi\in\SS_n$; we say $\pi_j$ is {\em
right to left maxima} of $\pi$ if $\pi_j>\pi_i$ for all $j<i$. The
number of right to left maxima of $\pi$ we denote by $rlm_\pi$.

\begin{corollary}\label{c1}
We have

{\rm(i)}
$\sum\limits_{n\geq0}\sum\limits_{\pi\in\EE_n(132)}x^ny^{rlm_\pi}=
\dfrac{1}{2}\left(\dfrac{1}{1-xyC(x)}+\dfrac{1+xy-x^2yC(x^2)}{1-2x^2yC(x^2)+x^2y^2C(x^2)}\right),$

{\rm(ii)}
$\sum\limits_{n\geq0}\sum\limits_{\pi\in\OO_n(132)}x^ny^{rlm_\pi}=
\dfrac{1}{2}\left(\dfrac{1}{1-xyC(x)}-\dfrac{1+xy-x^2yC(x^2)}{1-2x^2yC(x^2)+x^2y^2C(x^2)}\right).$

\end{corollary}
\begin{proof}
Using Theorem~\ref{thss} together with definitions we have
$\MM(x,1,1,\ldots)=1+xC(x^2)$. So, Theorem~\ref{thggg} yields
$$\MM(xy,y^{-1},y,\ldots)=\EE(xy,y^{-1},y,\ldots)-\OO(xy,y^{-1},y,\ldots)=\dfrac{1+xy-x^2yC(x^2)}{1-2x^2yC(x^2)+x^2y^2C(x^2)},$$
and
$$\F(xy,y^{-1},y,y^{-1},\ldots)=\EE(xy,y^{-1},y,y^{-1},\ldots)+\OO(xy,y^{-1},y,y^{-1},\ldots)=\frac{1}{1-xyC(x)}.$$
On the other hand, using \cite[Proposition~5]{BCS} we get
$$\begin{array}{l}
\sum\limits_{n\geq0}\sum\limits_{\pi\in\SS_n(132)}x^ny^{rlm_\pi}=\F(xy,y^{-1},y,y^{-1},\ldots),\\
\sum\limits_{n\geq0}\sum\limits_{\pi\in\EE_n(132)}x^ny^{rlm_\pi}=\EE(xy,y^{-1},y,y^{-1},\ldots),\\
\sum\limits_{n\geq0}\sum\limits_{\pi\in\OO_n(132)}x^ny^{rlm_\pi}=\OO(xy,y^{-1},y,y^{-1},\ldots).\\
\end{array}$$
By combining all these equations we get the desired result.
\end{proof}

An another application for Theorem~\ref{thggg} we get an explicit
expressions for the generating function
$\sum_{n\geq0}\sum_{\pi\in\EE_n(132)}x^ny^{12\ldots k(\pi)}$ for
given $k$ and $12\ldots k(\pi)$. The following result is true by
using Theorem~\ref{thggg}.

\begin{theorem}
Let $k\geq1$; we have

{\rm(i)}
$\MM_{[2k+1];0}(x)=1+xR_k(x^2)=\frac{U_k\sx+U_{k-1}\sx}{U_k\sx}$;

{\rm(ii)} $\MM_{[2k+1];1}(x)=\frac{x}{U_k^2\sx}$;

{\rm(iii)}
$\MM_{[2k+1];2}(x)=\frac{x^2}{U_k^2\sx}(xR_k(x^2)-1)=\frac{x^2(U_{k-1}\sx-U_k\sx)}{U_k^3\sx}$.
\end{theorem}

Therefore, by \cite[Theorem~4.1]{MV1} together with the above
theorem we get the number of even (or odd) permutations avoiding
$132$ and containing $[2k+1]$ exactly $r=0,1,2$. For example, for
$r=2$ we get (for $r=1$ see Corollary~\ref{ccra})

\begin{corollary}
Let $k\geq1$. Then

{\rm(i)} the generating function for the number $132$-avoiding
even permutations containing $12\ldots(2k+1)$ exactly twice is
given by
$$\frac{1}{2}\left(\frac{\sqrt{x}U_{k-1}\txx}{U_k^3\txx}+\frac{x^2(U_{k-1}\sx-U_k\sx)}{U_k^3\sx}\right),$$

{\rm(ii)} the generating function for the number $132$-avoiding
odd permutations containing $12\ldots(2k+1)$ exactly twice is
given by
$$\frac{1}{2}\left(\frac{\sqrt{x}U_{k-1}\txx}{U_k^3\txx}-\frac{x^2(U_{k-1}\sx-U_k\sx)}{U_k^3\sx}\right),$$
\end{corollary}

\subsection{Two restrictions}
The second of these directions is to consider more than
one additional restriction. For example, the following results is
true. Let $B_{\tau^1,\tau^2}(x)$ be the generating function for
the number of even permutations in $\EE_n(132,\tau^1,\tau^2)$.
Assume $\tau^1=12\ldots k$ and $\tau^2=2134\ldots k$, then we get
as follows.

\begin{theorem}\label{thja}
Let $V_{m}(x)=(1-xW_m(x))W_{m+1}(x)$ such that
$$W_m(x)=\frac{(1-x^2R_{m-2}(x)R_{m-3}(x))R_{m-1}(x)}{1-x^2R_{m-1}(x)R_{m-2}(x)}.$$
Then, for all $k\geq2$,
\begin{enumerate}
\item   The generating function $\EE_{12\ldots 2k,2134\ldots 2k}(x)$ is given by
            $$\dfrac{1}{2}\left(W_{2k}(x)+1+xR_k(x^2)\right).$$

\item   The generating function $\OO_{12\ldots 2k,2134\ldots 2k}(x)$ is given by
            $$\dfrac{1}{2}\left(W_{2k}(x)-1-xR_k(x^2)\right).$$

\item   The generating function $\EE_{12\ldots(2k+1),2134\ldots(2k+1)}(x)$ is given by
            $$\frac{1}{2}\left(\frac{2(1+xR_{k+1}(x^2))-V_{2k}(x)-xR_{k+1}(x^2)V_{2k}(-x)}{1+x^2R_{k+1}^2(x^2)}\cdot
                R_{k+1}(x^2)+W_{2k+1}(x)\right).$$

\item   The generating function $\OO_{12\ldots(2k+1),2134\ldots(2k+1)}(x)$ is given by
            $$\frac{1}{2}\left(W_{2k+1}(x)-\frac{2(1+xR_{k+1}(x^2))-V_{2k}(x)-xR_{k+1}(x^2)V_{2k}(-x)}{1+x^2R_{k+1}^2(x^2)}\cdot
                R_{k+1}(x^2)\right).$$
\end{enumerate}
\end{theorem}

Another example to consider the case of avoiding $\tau^1$ and
counting occurrences of $\tau^2$. For example, the following
result is true. Let $\EE_{\tau^1}^{\tau^2}(x,y)$
\resp{$\OO_{\tau^1}^{\tau^2}(x,y)$} be the generating function for
the number of even \resp{odd} permutations in $\SS_n(132,\tau^2)$
containing $\tau^1$ exactly $r$ times.

\begin{theorem}
Let $G_k(x,y)=\EE_{[k]}^{[k+1]}(x,y)-\OO_{[k]}^{[k+1]}(x,y)$. For
all $k\geq 1$,
$$G_k(x,y)=1+x\frac{D_{k-1}(x)-x^k+B_k(x)(1-y)+x^2(D_{k-3}(x)-x^{k-2})(1-y)^2}
{D_k(x)+E_k(x)(1-y)+x^2D_{k-2}(x)(1-y)^2},$$ where
$E_m(x)=((-1)^m-1)x^{m+1}$; $B_{2m}(x)=-x^{2m-1}$ and
$B_{2m+1}=-x^{2m}+2x^{2m+1}$; and
$$\begin{array}{l}
D_{2m}(x)=\frac{x^{2m+1}}{1-4x^2}(U_{2m+1}\sx-2xU_{2m}\sx-2x),\\
D_{2m+1}(x)=\frac{x^{2m+3}}{1-4x^2}(U_{2m+3}\sx-U_{2m+1}\sx-4x),
\end{array}$$
for all $m\geq0$.
\end{theorem}

{\bf Acknowledgments}. The final version of this paper was written
while the author was visiting University of Haifa, Israel in
January 2003. He thanks the HIACS Research Center and the Caesarea
Edmond Benjamin de Rothschild Foundation Institute for
Interdisciplinary Applications of Computer Science for financial
support, and professor Alek Vainshtein for his generosity.

\end{document}